\documentclass[11pt]{amsart}

\usepackage{a4wide}
\usepackage{euscript,amsmath,amssymb,amsbsy,mathabx}
\usepackage[arrow,matrix,curve]{xy}
\usepackage{array,longtable,graphicx,enumitem,url,multirow}

\usepackage{xcolor}
\usepackage[colorlinks=true,linkcolor=black,urlcolor=black,citecolor=black]{hyperref}

\usepackage{tikz,pgfplots}
\usetikzlibrary{patterns}
\pgfplotsset{compat=1.13}\usepgfplotslibrary{fillbetween}
\pgfdeclarelayer{ft}\pgfdeclarelayer{bg}\pgfsetlayers{bg,main,ft}

\newcounter{msct}[section]\renewcommand{\themsct}{\thesection.\arabic{msct}}
\newenvironment{m-theorem}{\vskip5pt\refstepcounter{msct}\trivlist \itemindent 0pt%
\item[\hskip\labelsep\bf Theorem~\themsct]\it\ignorespaces}{\endtrivlist\vskip3pt}
\newenvironment{m-proposition}{\vskip5pt\refstepcounter{msct}\trivlist \itemindent0pt%
\item[\hskip\labelsep\bf Proposition~\themsct]\it\ignorespaces}{\endtrivlist\vskip3pt}
\newenvironment{m-corollary}{\vskip5pt\refstepcounter{msct}\trivlist \itemindent 0pt%
\item[\hskip\labelsep\bf Corollary~\themsct]\it\ignorespaces}{\endtrivlist\vskip3pt}
\newenvironment{m-lemma}{\vskip5pt\refstepcounter{msct}\trivlist \itemindent 0pt%
\item[\hskip\labelsep\bf Lemma~\themsct]\it\ignorespaces}{\endtrivlist\vskip3pt}
\newenvironment{m-definition}{\vskip5pt\refstepcounter{msct}\trivlist \itemindent0pt%
\item[\hskip\labelsep\bf Definition~\themsct]\ignorespaces}{\endtrivlist\vskip5pt}
\newenvironment{m-notation}{\vskip5pt\refstepcounter{msct}\trivlist \itemindent0pt%
\item[\hskip\labelsep\bf Notation~\themsct]\ignorespaces}{\endtrivlist\vskip5pt}
\newenvironment{m-example}{\vskip5pt\refstepcounter{msct}\trivlist \itemindent0pt%
\item[\hskip\labelsep\bf Example~\themsct]\ignorespaces}{\endtrivlist\vskip5pt}
\newenvironment{m-remark}{\vskip5pt\refstepcounter{msct}\trivlist \itemindent0pt%
\item[\hskip\labelsep\bf Remark~\themsct]\ignorespaces}{\endtrivlist\vskip5pt}
\newenvironment{m-question}{\vskip3pt\refstepcounter{msct}\trivlist \itemindent0pt%
\item[\hskip\labelsep\bf Question.]\ignorespaces}{\endtrivlist\vskip5pt}
\newenvironment{thm-nono}[1]{\vskip5pt\trivlist \itemindent 0pt %
\item[\hskip\labelsep\bf Theorem~{\rm\mbox{#1}}]\it\ignorespaces}{\endtrivlist\vskip5pt}
\newenvironment{prop-nono}[1]{\vskip5pt\trivlist \itemindent 0pt %
\item[\hskip\labelsep\bf Proposition~{\rm\mbox{#1}}]\it\ignorespaces}{\endtrivlist\vskip5pt}
\newenvironment{lm-nono}[1]{\vskip5pt\trivlist \itemindent0pt%
\item[\hskip\labelsep\bf Lemma~{\rm\mbox{#1}}]\it\ignorespaces}{\endtrivlist\vskip5pt}
\newenvironment{conj-nono}[1]{\vskip5pt\trivlist \itemindent0pt%
\item[\hskip\labelsep\bf Conjecture~{\rm\mbox{#1}}]\it\ignorespaces}{\endtrivlist\vskip5pt}
\newenvironment{def-nono}[1]{\vskip5pt\trivlist \itemindent0pt%
\item[\hskip\labelsep\bf Definition~{\rm\mbox{#1}}]\ignorespaces}{\endtrivlist\vskip5pt}
\newenvironment{m-thank}{\vskip5pt\trivlist \itemindent0pt%
\item[\hskip\labelsep\it Acknowledgments]\ignorespaces}{\endtrivlist\vskip5pt}
\newenvironment{m-proof}{\vskip2pt\trivlist \itemindent0pt%
\item[\hskip\labelsep\it Proof.]\ignorespaces}{\hfill$\Box$\endtrivlist\vskip5pt}%
\newenvironment{m-asmp}{\vskip5pt\refstepcounter{msct}\trivlist \itemindent0pt%
\item[\hskip\labelsep\bf Assumption~\themsct]\ignorespaces}{\hfill\endtrivlist\vskip5pt}%

\newcounter{meqn}[section]\renewcommand{\themeqn}{\thesection.\arabic{meqn}}
\newenvironment{m-eqn}[1]{\vskip5pt\refstepcounter{meqn}\trivlist\itemindent0pt%
\item[]\ignorespaces\hfill$\displaystyle #1$\hfill\hbox{\rm(\themeqn)}}{\endtrivlist\vskip5pt}

\newcommand{\bibauth}[2]{\textrm{{#1}~{#2}}}
\newcommand{\bibtitl}[1]{\textit{#1}.}
\newcommand{\bibjnyp}[4]{\textrm{#1} \textbf{#2} (#3), #4.}
\newcommand{\bibinbook}[4]{In: \textrm{#1}\textrm{, #2}\textrm{, #3}\textrm{, #4}.}
\newcommand{\bibbook}[4]{\textit{#1}. {#2} {#3}, {#4}.}

\numberwithin{equation}{section}\numberwithin{figure}{section}

\let\euf\EuScript 
\let\cal\mathcal
\let\mbb\mathbb
\let\mfrak\mathfrak 

\newcommand\ouset[3]{{\overset{#2}{\underset{#1}#3}\,}}

\let\dta\delta
\let\disp\displaystyle 
\newcommand{\eI}{{\euf I}}
\newcommand{\cd}{\mathop{\rm cd}\nolimits}
\renewcommand{\dim}{\mathop{\rm dim}\nolimits}
\newcommand{\codim}{\mathop{\rm codim}\nolimits}
\newcommand{\Gr}{\text{\rm Gr}}
\newcommand{\oGr}{\text{\rm o-Gr}}
\newcommand{\sGr}{\text{\rm sp-Gr}}
\let\ges\geqslant 
\let\les\leqslant 
\let\nit\noindent 
\newcommand{\lcit}{{\textit{loc.\,cit.}}}

\newcommand{\kk}{{\Bbbk}}
\newcommand{\eE}{{\euf E}}
\newcommand{\eF}{{\euf F}}
\newcommand{\eG}{{\euf G}}
\newcommand{\eL}{{\euf L}}
\newcommand{\bbM}{{\mbb L}}
\newcommand{\eM}{{\euf L}}
\newcommand{\main}{{\mathrm{main}}}
\newcommand{\mn}{{\rm min}}
\newcommand{\mx}{{\rm max}}
\newcommand{\eN}{{\euf N}}
\newcommand{\eO}{{\euf O}}
\newcommand{\eQ}{{\euf Q}}
\newcommand{\pnt}{{pt}}
\newcommand{\msp}{{\mfrak L}}
\newcommand{\spp}{{\mfrak{l}}}

\newcommand{\eT}{{\euf T}}
\newcommand\Hom{\mathop{\rm Hom}\nolimits}
\newcommand\Pic{\mathop{\rm Pic}\nolimits}
\newcommand{\Sym}{\mathop{\rm Sym}\nolimits}
\newcommand{\ut}{{\underline{t}}}
\let\ovl\overline 
\let\sm\setminus 
\let\srel\stackrel 
\let\tld\tilde 
\let\unbar\underbar

\begin{document}

\title[Criteria for complete intersections]{Criteria for complete intersections}
\author{Mihai Halic}
\email{mihai.halic@gmail.com}

\begin{abstract}
We obtain criteria for detecting complete intersections in projective varieties. Motivated by a conjecture of Hartshorne concerning subvarieties of projective spaces, we investigate situations when two-co\-dimen\-sional smooth subvarieties of rational homogeneous varieties are complete intersections.
\end{abstract}

\maketitle

\section*{Introduction}

Hartshorne~\cite{hart-conj} conjectured that smooth subvarieties of projective spaces of sufficiently low codimension are necessarily complete intersections. He also stated that this is indeed the case as soon as the degree of the subvariety is small compared to the dimension of the ambient projective space.

The goal of this article is to investigate this problem in a more general setting. We ask whether a similar statement holds for smooth subvarieties $X$ of other ambient varieties $V$ as well. From the beginning, one distinguishes two cases: 
\begin{itemize}[leftmargin=4ex]
\item 
One assumes that the normal bundle of $X$ decomposes into a direct sum of line bundles. This condition is clearly satisfied for complete intersections.
\item[or]
\item 
There are no \textit{a priori} hypotheses. 
\end{itemize}
In the first case, the property of being a complete intersection is reduced to verifying the splitting of the normal bundle of the subvariety; this feature is, at least theoretically, much easier to verify. Faltings~\cite{falt-krit} showed that a smooth subvariety of the projective space $\mbb P^n$, of codimension at most $n/2$, is a complete intersection as soon as its normal bundle splits. By analysing his proof, one observes that at its base there is a Hauptlemma---main lemma---which is valid in great generality. However, its `implementation' is strongly adapted to subvarieties of projective spaces. This may explain why, to our knowledge, the arguments have not been extended to a wider framework. In the first part of the article we attempt bridging this difference. The following statement conveys the flavour of our results. 

\begin{thm-nono}{(cf.~\ref{thm:homog},~\ref{thm:generators})} 
Let $G$ be a simple, linear algebraic group of rank $\ell\ges6$ and $P$ a maximal parabolic subgroup. 

Let $X$ be a smooth subvariety of the rational homogeneous variety $G/P$, of codimension $2\les\dta\les(\ell+1)/3$. Then $X$ is a complete intersection as soon as either one of the following two conditions is satisfied: 
\begin{itemize}[leftmargin=4ex]
\item 
Its normal bundle $\eN_{X/V}$ splits into a direct sum of line bundles.
\item 
The ideal defining $X$ is generated by at most $(\ell+1)/2-\dta$ equations. 
\item[] (By equations, we mean sections in various $\eO_V(d),\;d>0.$)
\end{itemize}
\end{thm-nono}
In down-to-earth terms, the splitting of the normal bundle of the subvariety or the generation of its defining ideal by few enough polynomials ensure the complete intersection property. Let us remark that we prove this kind result for a much wider class of situations (cf.~Theorem~\ref{thm:s-split}); the formulation above is chosen to avoid additional terminology. 

The second part is devoted to the case where there are no assumptions on the splitting type of the normal bundle of the subvariety; we restrict ourselves to the two-codimensional case. For subvarieties of arbitrary codimension in projective spaces, positive results towards Hartshorne's conjecture have been obtained by Barth-Van de Ven~\cite{bart-icm} and Bertram-Ein-Lazarsfeld~\cite{bel}. In the two-codimensional case, Barth-Van de Ven~\cite{bart+vdv}, Ran~\cite{ran}, Holme-Schneider~\cite{holm+sch} determined effective upper bounds for the degree of the subvariety ensuring the complete intersection property. Also, the conjecture has been proved by Ionescu-Russo~\cite{ion+rus} for subvarieties cut out by quadratic equations. 

\begin{thm-nono}{(cf.~\ref{thm:ran})} 
Let $G/P$ be as before, with $\ell\ges6$, and denote $m:=\deg_L(\eT_V)-3$, where $L$ is the $1$-dimensional Schubert line. Suppose $X$ is smooth, $2$-codimensional, such that $[X]=d_X\cdot[\eO_{G/P}(1)]^2,\;\det(\eN_{X/V})=\eO_X(n_X).$ If either 
$$
d_X\les m(n_X-m)\quad\text{or}\quad\sqrt{d_X}\les n_X/2\les m,
$$
the subvariety $X\subset G/P$ is a complete intersection. 
\end{thm-nono}

The statement is weaker than before, when the splitting of the normal bundle was known in advance, but it faithfully parallels the results of Ran and Schneider \emph{for subvarieties of projective spaces}. In fact, along the way, we show (cf. Theorem~\ref{thm:d>}) that subvarieties with stable normal bundle have large degree. Our bound (of order four in $\ell$) is better than that of Schneider (quadratic in $\ell$), even for projective spaces. This fact, combined with the theorem above yields a positive answer to Hartshorne's question for subvarieties of low degree. 
\begin{thm-nono}{(cf.~\ref{thm:hart})} 
Let $G/P$ be as before, $\ell\ges11$, and $m:=\deg_L(\eT_V)-3$. Let $X$ be smooth, $2$-codimensional, such that 
$$[X]=d_X\cdot[\eO_{G/P}(1)]^2,\qquad d_X\les m^2.$$  
Then $X$ is a complete intersection in $V$.
For $\ell=7,\dots,10$, the same holds for $d_X\les3m^2/10$.
\end{thm-nono}

To the author's knowledge, the \emph{only one reference} which investigates the complete intersection property for subvarieties of homogeneous varieties, other than projective spaces, is Barth-Van de Ven~\cite{bart+vdv2} dealing with $2$-codimensional subvarieties of Grassmannians (of very low degree). 

We believe that the interest in our work (with or without assumptions on the normal bundle) consists in placing Hartshorne's question in a much wider context.


\section{Subvarieties with split normal bundle}\label{sct:proof}

\begin{m-notation}\label{not:XV} 
We work over an algebraically closed field $\kk$ of characteristic zero. Throughout the article, $X$ is a smooth subvariety of an ambient smooth projective variety $V$, of codimension $2\les\dta\les\dim V/2$, which is defined by the sheaf of ideals $\eI_X\subset\eO_V$; its co-normal bundle  $\eN_{X/V}^\vee:=\eI_X/\eI_X^2$ is locally free, of rank $\dta$. For a subvariety $Y$ of $V$ containing $X$, we denote $\eI_{X\subset Y}:=\eI_X/\eI_Y$. 
\end{m-notation} 
In this section, we consider the following issue: 

\medskip Assuming that  $\eN_{X/V}^\vee$ splits ---it decomposes into a direct sum of line bundles---, can one deduce that $X$ a complete intersection of $\dta$ hypersurfaces in $V$?\medskip  

Faltings' work~\cite{falt-krit} corresponds to the case where $V$ is a projective space. 

Our investigations naturally involve partially ample subvarieties of projective varieties (investigated by the author in~\cite{hlc}). For the reader's comfort, we briefly recalled these notions in the Appendix~\ref{sct:partialpos}. 

We assume henceforth, throughout this section, that the following conditions are satisfied.

\begin{m-asmp}\label{asmp:XV}
\begin{enumerate}[leftmargin=4ex]
\item 
The conormal bundle of $X$ in $V$ splits into a direct sum of line bundles
\begin{m-eqn}{
\eI_X/\eI_X^2=\ouset{j=1}{\dta}{\bigoplus}\eL_j^{-1}\otimes \eO_X,\quad \eL_j\in\Pic(V).
}\label{eq:N-split}
\end{m-eqn}
We denote by 
\begin{m-eqn}{
x_j\in\Gamma\big(X,(\eI_X/\eI_X^2)\otimes \eL_j\big),\;1\les j\les\dta,
}\label{eq:x}
\end{m-eqn}
the section corresponding to the direct summand $\eO_X$.\smallskip  
\item 
$\,\eL_{j-1}^{-1}\eL_j$ is semi-ample, for all $j=1,\dots,\dta$,. (We agree that $\eL_0:=\eO_V$.) 
Note that, in this case, $\eL_j$ is semi-ample too. For short, we say that $\eL_1,\dots,\eL_\dta$ are \emph{ordered}. \smallskip 
\item 
The subvariety $X\subset V$ is $q_X$-ample, $q_X\les\dim V-2\dta=\dim X-\dta$ (so $\dim V\ges4$).
\item[] 
Together with the previous assumption, this implies that $\eL_j\otimes\eO_X$ is $q_X$-positive, for all $j$. (See Proposition~\ref{prop:XV} and Definition~\ref{def:part-pos}.)
\end{enumerate}
\end{m-asmp}

\begin{m-remark} 
\begin{enumerate}[leftmargin=4ex]
\item 
We imposed that the line bundles appearing in the splitting of the conormal bundle (over $X$) actually come from $V$. This is not very far from the full generality, where they would be defined only over $X$: for $\dta\ges3$ and arbitrary $V$, the partial amplitude of $X$ implies that $\Pic(V)\otimes\mbb Q\to\Pic(X)\otimes\mbb Q$ is an isomorphism (cf.~\cite[Corollary 1.12]{hlc}). 
Actually the isomorphism holds with $\mbb Z$-coefficients in many situations, especially for subvarieties of rational homogeneous varieties and for zero loci of sections in vector bundles (cf. Barth-Lefschetz-type theorems~\cite{soms,soms-vdv}). 
\item 
For rational homogeneous varieties $V=G/P$, where $P\subset G$ is a maximal parabolic subgroup (so $\Pic(V)\cong\mbb Z$), the line bundles $\eL_j$ are positive multiples of $\eO_{V}(1)$, hence $q_j=0$ for all $j$. Consequently, our partial positivity assumption is a weakening of this situation. 
\item
The semi-ampleness of the successive differences means that $\eL_1,\dots,\eL_\dta$ are ordered `increasingly'; the condition is automatically satisfied whenever the Picard group is cyclic.
\item 
The most important information given by the partial amplitude of $X\subset V$ is the upper bound on the cohomological dimension of its complement (see Proposition~\ref{prop:XV}). It ensures that certain intersections are non-trivial, substituting the knowledge of the cohomology ring of $V$. 
\end{enumerate}
\end{m-remark}

\begin{m-lemma}\label{lm:inters}
Let $V,X$ be as above.  Suppose $Z\subset V$ is an arbitrary closed, irreducible subscheme of codimension at most $\dta$. Then $Z$ intersects $X$ non-trivially.
\end{m-lemma}

\begin{m-proof}
The partial amplitude of $X\subset V$ yields: 
$\cd(V\sm X)\les q_X+\dta-1\les\dim X-1.$ 
Hence $Z$, which is complete with $\dim Z\ges\dim X$, can't be contained in the complement of $X$. 
\end{m-proof}

\begin{m-lemma}\label{lm:surj}
Let $V,X$ be as above. Let $Y$ be an irreducible, normal, lci subvariety of $V$ which (strictly) contains $X$, is smooth along it, and such that  
$$
\eN^\vee_{X/Y}=\eI_{X\subset Y}/\eI_{X\subset Y}^2=\ouset{j=d}{\dta}{\bigoplus}\eL_j^{-1}\otimes \eO_X 
$$
for some $d\les\dta$. Then the following statements hold.
\begin{enumerate}[leftmargin=4ex]
\item 
There is a unique $y\in\Gamma(Y,\eI_{X\subset Y}\otimes \eL_{d})$ which is mapped to $x_d\in\Gamma(X,\eN_{X/V}^\vee\otimes \eL_{d})$ under the natural homomorphism $\Gamma(Y,\eL_{d})\to\Gamma(X,(\eO_Y/\eI_{X\subset Y}^{2})\otimes \eL_{d}).$
\item 
Let $Z\subset Y$ be the vanishing locus of the section $y$. Then $Z$ is still an irreducible, normal, lci variety containing $X$, and is smooth along it; moreover, we have 
$$
\eN^\vee_{X/Z}=\eI_{X\subset Z}/\eI_{X\subset Z}^2=\ouset{j=d+1}{\dta}{\bigoplus}\eL_j^{-1}\otimes \eO_X.
$$
\end{enumerate}
\end{m-lemma}

\begin{m-proof}
(i) Let $\hat Y_X$ be the formal completion of $Y$ along $X$. The statement is deduced in two steps. 
\begin{itemize}[leftmargin=4ex]
\item 
We claim that the homomorphisms
$$
\Gamma\big(Y,(\eO_Y/\eI_{X\subset Y}^{r+1})\otimes \eL_{d}\big)\to \Gamma\big(Y,(\eO_Y/\eI_{X\subset Y}^{r})\otimes \eL_{d}\big),\quad r\ges 2,
$$
are isomorphisms. Indeed, we have the exact sequence 
$$
0\to \Sym^r\biggl(\frac{\eI_{X\subset Y}}{\eI^2_{X\subset Y}}\biggr)\otimes \eL_{d}
\to\frac{\eO_Y}{\eI^{r+1}_{X\subset Y}}\otimes \eL_{d}
\to \frac{\eO_Y}{\eI^{r}_{X\subset Y}}\otimes \eL_{d}\to 0.
$$
The assumption on the conormal bundle implies that the left-hand side is a direct sum of line bundles $M^{-1}$ of the form 
$$
M=\big(\eL_{j_1}\otimes\dots\otimes \eL_{j_r}\otimes \eL_{d}^{-1}\big)\otimes\eO_X,\quad j_1,\dots,j_r\ges d.
$$ 
Since $\eL_{j_1}\eL_{d}^{-1}$ is semi-ample and $r\ges2$, each such $M$ is $q_{j_2}$-positive, with $q_{j_2}\les\dim X-2$. The Grauert-Riemenschneider theorem yields the vanishing of the cohomology groups of $M^{-1}$ in degrees zero and one. 

\item 
The second claim is that $\Gamma(Y,\eL_{d})\to \Gamma(\hat Y_X, \hat\eL_d)$ is an isomorphism, $\hat\eL_d:=\eL_{d}\otimes\eO_{\hat Y_X}$.
\item[] It is enough to show that, for all sufficiently large integers $r$, it holds 
$$
H^0(Y,\eI_{X\subset Y}^r\otimes\eL_{d})=H^1(Y,\eI_{X\subset Y}^r\otimes\eL_{d})=0.
$$ 
By using the definition, we see that $X\subset Y$ is $(q_X+d-1)$-ample and $q_X+d-1\les\dim X-1$. Since $Y$ is a Gorenstein variety (it is lci) and is smooth along $X$, the blow-up $\tld Y$ of $X$ is still Gorenstein, so it admits a dualising sheaf $\omega_{\tld Y}$; let $E_X$ be the exceptional divisor. Then we have, $\text{for}\;t=0,1\;\text{and}\;r\gg0,$ 
$$
\begin{array}{r}
H^t(Y,\eI_{X\subset Y}^r\otimes\eL_{d})=H^{\dim Y-t}(\tilde Y,\eO_{\tld Y}(-rE_X)\otimes \eL_{d})
\\[1.5ex] 
=H^{\dim Y-t}(\tilde Y,\omega_{\tld Y}\otimes\eO_{\tld Y}(rE_X)\otimes \eL_{d})=0.
\end{array}
$$
\end{itemize}
Overall, $\Gamma(Y,\eL_{d})\to\Gamma(X,(\eO_Y/\eI_{X\subset Y}^2)\otimes \eL_{d})$ is an isomorphism.\smallskip

\nit(ii) The scheme $Z$ contains $X$, by construction. The section $x_d$ can be viewed as the differential of $y$ along its vanishing locus, so the differential criterion implies that $Z$ is smooth along $X$. Thus there is a unique irreducible component $Z_o\subset Z$ containing $X$. 

Since $Z$ is the zero locus of a section in a line bundle on $Y$, all its components are hypersurfaces in $Y$. Let $Z'$ be an arbitrary component, so $\codim_VZ'=d$. Lemma~\ref{lm:inters} implies that $Z'\cap X\neq\emptyset$, so $Z'\cap Z_o\neq\emptyset$; the previous discussion yields $Z=Z_o$.

It remains to show that $Z$ is normal; it's enough to prove that the singular locus $Z^{sing}\subset Z$ has codimension at least two. If  $\dim Y=\dim X+1$, then $Z_o=X$ is smooth. Now suppose $\dim Y\ges\dim X+2$. Then $Z^{sing}$ is disjoint of $X$, so we have 
$$
\dim Z^{sing}\les\cd(Y\sm X)\les\dim X-1\les\dim Z-2.
$$ 
The last statement follows from the exact sequence of co-normal bundles. 
\end{m-proof}

By inductively applying the previous lemma, we obtain the following statement which holds in great generality. 

\begin{m-proposition}\label{prop:flag}
Suppose $V,X$ are as in~\ref{asmp:XV}. Then there is a flag of irreducible, normal, lci (hence Cohen-Macaulay) subvarieties 
$$
V=Y_0\supset Y_1\supset\dots\supset Y_{\dta-1}\supset Y_{\dta}=X.
$$
Each $Y_j\subset Y_{j-1}$ is the vanishing locus of a section $y_j\in\Gamma(Y_{j-1},\eL_j)$, for $j=1,\dots,d$, which extends $x_j\in\Gamma(X,\eN^\vee_{X/V}\otimes\eL_j)$. 
\end{m-proposition}

We continue our investigation, because we are interested whether $X$ is a complete intersection. The issue lies at Lemma~\ref{lm:surj}, which doesn't extend the sections ${\{x_j\}}_j$ to the whole ambient space $V$. Let us show that this is the only obstruction to the desired conclusion.

\begin{m-lemma}\label{lm:XV}
Let $d\les\dta$ and $v_1\in\Gamma(V,\eI_X\otimes \eL_1),\dots,v_d\in\Gamma(V,\eI_X\otimes \eL_d)$ which induce the sections $x_1,\dots,x_d$, respectively; let $Y_d$ be their common (scheme-theoretic) vanishing locus. Then $Y_d$ is an irreducible, normal, $d$-codimensional complete intersection in $V$.

Suppose the ground field is $\mbb C$. The same conclusion holds if, instead of the Assumption~\ref{asmp:XV}(iii), $X$ satisfies the following condition. 
\\ Its fundamental class is Poincar\'e-dual to $d_X\cdot\chi^\dta,\;d_X>0,$ and $\chi\in H^2(V;\mbb Z)$ is a K\"ahler class. (More generally, $[X]$ belongs to the positive cone generated by K\"ahler classes.)
\end{m-lemma}

\begin{m-proof}
By their very definition, the divisors $D_j:=\{v_j=0\}$ contain $X$, so $Y_d$ contains $X$ too. Moreover, the differential criterion shows that $D_1,\dots,D_d$ intersect transversally along $X$, therefore $Y_d$ is smooth in a neighbourhood of $X$. We deduce that there is a unique irreducible component $Y\subset Y_d$ containing $X$ which has codimension $d$ in $V$. 

We are going to show that $Y_d$ is irreducible. Take some other component $Y'$ of $Y_d$; its dimension is at least $\dim V-d$ and Lemma~\ref{lm:inters} shows that $Y'\cap X\neq\emptyset$. The smoothness of $Y_d$ along $X$ implies that $\codim_V(Y')=d$ and $Y'=Y$. Thus $Y_d$ is an irreducible, complete intersection subvariety of $V$. 
For $d=\dta$, we have $Y=X$ is smooth. For $d\les\dta-1$, the upper bound on the cohomological dimension of $V\sm X$ forces the singular locus $Y_d^{sing}$, which is disjoint of $X$, to have $\dim Y_d^{sing}\les\dim X-1\les\dim Y_d-2$. It follows that $Y_d$ is a normal variety. 

Now we turn to the second statement. The previous proof applies \textit{ad litteram} as long as Lemma~\ref{lm:inters} is valid: $X$ intersects any closed subscheme $Z$ of codimension at most $\dta$. We may assume that $Z$ is reduced, irreducible, $\dta$ codimensional, so its Poincar\'e dual belongs to $H^{2\dta}(V;\mbb R)$. Then $[Z]\cap(\chi^\dta\cup\chi^{\dim V-2\dta})$ is the volume of $Z$ with respect to the K\"ahler form, a strictly positive number. This proves that $[Z]\cap[X]\neq0$. 
\end{m-proof}

Extending the sections $y_j\in\Gamma(Y_{j-1},\eL_{j})$ to $V$ require stronger cohomology vanishing properties for line bundles on $V$; the assumptions which will be made in the sequel are imposed by the analysis of the general case. Recall that $v_1:=y_1$ is already defined over $V$. Suppose that $v_1,\dots,v_{d-1},\;d\les\dta,$ as above are already constructed, so we have the Koszul resolution 
$$
0{\to}\ouset{}{d-1}{\bigwedge}\eN_{d-1}^\vee{\to}\ouset{}{d-2}{\bigwedge}\eN_{d-1}^\vee{\to}\dots{\to}\eN_{d-1}^\vee{\to}\eO_X{\to}\eO_{Y_{d-1}}{\to} 0,\quad \eN_{d-1}:=\eL_1\oplus\dots\oplus\eL_{d-1}. 
$$
For lifting $y_d$, it suffices to ensure that $H^t\Big(V,\ouset{}{t}{\bigwedge}\eN_{d-1}^\vee\otimes\eL_{d}\Big)=0,\; 1\les t\les d-1.$ As $d$ varies, within the brackets appear direct sums of the form $\eL_{j_1}^{-1}\cdot\dots\cdot\eL_{j_t}^{-1}\eL_{d},$ so one should have: 
\begin{m-eqn}{
\forall\; j_1<\dots<j_t<d\les\dta,\quad 
\begin{array}{l}
H^{t}\Big(V,\eL_{j_1}\cdot\dots\cdot\eL_{j_t}\eL_{d}^{-1}\Big)=0,\; \text{or equivalently}, 
\\[2ex]
H^{\dim V-t}\Big(V,\big(\omega_V^{-1}\otimes \eL_{j_1}^{-1}\cdot\dots\cdot\eL_{j_t}^{-1}\eL_{d}\big)^{-1}\Big)=0.
\end{array}
}\label{eq:H=0}
\end{m-eqn}
Now we distinguish two alternatives to achieve this matter: 
\begin{itemize}[leftmargin=4ex]
\item 
\unbar{either} by imposing more conditions on the line bundles $\eL_j$; 
\item 
\unbar{or} by restricting the type of varieties $V$ and keeping the line bundles $\eL_j$ arbitrary. 
\end{itemize}
Faltings' work concerns the projective space, so it enters into the latter category. To justify the following definition, observe that $\mbb P^n$ has the special feature that the line bundles on it have no intermediate cohomology.

\begin{m-definition}\label{def:s-split}
We say that the variety $V$ is \emph{$s$-split}, for an integer $0<s<\dim V$, if it satisfies the following property:
$$
H^t(V,\eL)=0,\quad\forall 1\les t\les s,\;\forall \eL\in\Pic(V).
$$ 
\end{m-definition}

Any $n$-dimensional Fano variety $V$ with cyclic Picard group is $(n-1)$-split. This property holds, in particular, for rational homogeneous varieties $G/P$, where $G$ is simple and $P$ is a maximal parabolic subgroup.  More generally, if $V$ is $s$-split and $W$ is a complete intersection of $\dta$ ample divisors, $\dim W\ges 3$, then $W$ is $(s-\dta)$-split. 

\begin{m-theorem}\label{thm:s-split} 
Let $V,X$ be as in Assumption~\ref{asmp:XV}. Suppose moreover that one of the following two properties is satisfied:
\begin{enumerate}[leftmargin=4ex]
\item 
The ambient variety $V$ is $s$-split and $X$ has codimension $\dta\les s-1$ in $V$.
\item 
The line bundle $\omega_V^{-1}\eL_1^{-1}\cdot\dots\cdot\eL_{\dta-1}^{-1}$ is ample; in particular, $V$ is a Fano variety.
\item[] 
(Loosely speaking, $X$ is a sufficiently low degree subvariety of $V$.)
\end{enumerate}
Then $X$ is a complete intersection in $V$. 
\end{m-theorem}

\begin{m-proof} In both cases we claim that the vanishing~\eqref{eq:H=0} is satisfied, so the sections $y_1,\dots,y_\dta$ extend to the ambient space $V$ and Lemma~\ref{lm:XV} applies. In the first case there is nothing to prove.

For the second one, we apply the Kodaira vanishing theorem. The successive differences $\eL_{j-1}^{-1}\eL_j$ are semi-ample, so it's enough to have the ampleness of  $\omega_V^{-1}\eL_{j_1}^{-1}\cdot\dots\cdot\eL_{j_{t-1}}^{-1}$; this follows from the hypothesis, since $\eL_j,\,j\not\in\{j_1,\dots,j_{t-1}\},$ are semi-ample. 
\end{m-proof}


\section{Applications}\label{sct:expl}

We claimed in the Introduction that, by considering partially ample subvarieties, the range of applications is substantially extended. Now we justify this claim and illustrate our result in several situations. 

\subsection*{Two-codimensional subvarieties}\label{ssct:2cod}

The interest in this case stems from Hartshorne's conjecture saying that the complete intersections are the only such subvarieties of $\mbb P^n,\,n\ges6$. However, the same issue can be raised more generally, so one may wonder what information can be extracted from the present work.  

Theorem~\ref{thm:s-split} yields a flag $V\supset Y_1\supset Y_2=X$ of length two. The obstruction for obtaining a complete intersection is the lifting of the section $y_2\in\Gamma(Y_1,\eL_2)$ to the ambient space $V$, so it is sufficient to ensure the vanishing of $H^1(V,\eL_1^{-1}\eL_2)$ whenever $\eL_1^{-1}\eL_2$ is semi-ample; one obtains further criteria, specific to this case.

\begin{m-corollary}\label{cor:codim2}
Suppose $V$ satisfies either one of the following conditions: 
\begin{itemize}[leftmargin=4ex]
\item it is $1$-split;
\item it is a Fano variety;
\item it is Frobenius split, compatible with an ample divisor (\textit{e.g.} toric, spherical variety, cf.~\cite{br+in}). 
\end{itemize}
Let $X\subset V$ be a smooth, $2$-codimensional, $(\dim X-2)$-ample subvariety, such that  
$$
\eN_{X/V}=(\eL_1\oplus\eL_2)\otimes\eO_X,\quad \eL_1,\eL_2\in\Pic(V),\;\text{with $\eL_1^{-1}\eL_2$ semi-ample.}
$$
Then $X$ is a complete intersection. 
\end{m-corollary}


In this statement, the partial ampleness of $X$ requires that the cohomological dimension of $V\sm X$ is at most $\dim V-3$. This may be difficult to verify in general, so it's desirable avoiding such assumptions. 

\begin{m-theorem}\label{thm:g3}
Let $V$ be a Fano variety with $\Pic(V)=\mbb Z\cdot\eO_V(1)$. Suppose that $X$ has the following properties:
\begin{itemize}[leftmargin=4ex]
\item 
$\eN_{X/V}=(\,\eO_V(a_1)\oplus\eO_V(a_2)\,)\otimes\eO_X,\quad0<a_1\les a_2;$
\item 
$[X]$ is Poincar\'e-dual to a multiple of $[\eO_V(1)]^2$;
\item 
$X$ is movable that is, its deformations cover an open subset of $V$.
\end{itemize}
Then $X$ is a complete intersection in $V$.
\end{m-theorem}

\begin{m-proof}
Let $x_1,x_2$ be defined as in~\eqref{eq:x}; we lift them from $X$ to $V$ in two steps, as in the proof of Lemma~\ref{lm:surj}. 

For $x_1$, the first step is identical: there is a lifting $\hat x_1$ to $\hat V_X$. The change comes to prove that $\Gamma(V,\eO_V(a_1))\to \Gamma(\hat V_X,\hat \eO_V(a_1))$ is an isomorphism. After dividing by a section in $\eO_V(a_1)$, $\hat x_1\in\Gamma(\widehat{\eO_V}(a_1))$ can be viewed as a formal rational function along $X$. 

Let $K(\hat V_X)$ be the field of formal rational functions; $\hat x_1$ lifts to $v_1\in\Gamma(V,\eO_V(a_1))$ if one has the equality $K(V)=K(\hat V_X)$, that is $X$ is G3 in $V$. For deducing this property, we apply~\cite[Theorem 8.2]{hlc-q} (cf. also proof of~\cite[Theorem 13.4(ii)]{bad}): 
\begin{itemize}[leftmargin=4ex]
\item the normal bundle of $X$ is ample, so it is G2 in $V$ (cf.~\cite{hlc});
\item $X$ is movable and intersects all the divisors of $V$;
\item $V$ is algebraically simply connected because it is Fano (hence rationally connected).
\end{itemize}

Now we turn to the section $x_2$. It suffices to lift it to $Y$, the zero set of the section $v_1$; recall that it is smooth along $X$. Once more, $x_2$ lifts to $\hat Y_X$ and it is enough to show that $X$ is G3 in $Y$. The normal bundle $\eN_{X/Y}=\eO_V(a_2)\otimes\eO_X$ is ample, so the G3 property follows from~\cite[Corollary 9.27]{bad}. 
\end{m-proof}


\subsection{Subvarieties of homogeneous varieties}\label{ssct:homog}

Here we work over the field of complex numbers. Suppose $V=G/P$, where $G$ is a semi-simple rational algebraic group and $P$ is a parabolic subgroup, and let  
$$
\ell:=\mn\{\;\text{rank}(H)\mid H\; \text{is a simple factor of}\;G\;\}.
$$
Note that $V$ decomposes into the product $\underset{H}{\prod}\,(H/P_H)$, where $H$ runs over the simple factors of $G$, and $P_H\subset H$ are parabolic subgroups.

A crucial ingredient is Faltings' upper bound for the cohomological dimension of the complement of $X$ and for the amplitude of the tangent bundle of $V$: 
\begin{m-eqn}{
\cd(V\sm X)\les\dim V-\ell+2\dta-2,\qquad\eT_V\;\text{is $(\dim V-\ell)$-ample.}
}\label{eq:faltings}\end{m-eqn}
Thus the direct summands of $\eN_{X/V}$ are $(\dim V-\ell)$-ample, too, and $X$ is a $q_X$-ample subvariety, with $q_X=\dim X-(\ell-2\dta+1)$, see Proposition~\ref{prop:XV}. 
Now we impose the conditions~\ref{asmp:XV} to obtain the upper bound on the codimension of $X$ for which Theorem~\ref{thm:s-split} applies:
$$
\begin{array}{rcl}
q_X\les\dim X-\dta&\quad\Rightarrow\quad&\dta\les(\ell+1)/3, 
\\[0ex]
\dta\ges2 &\quad\Rightarrow\quad& \ell\ges 5.
\end{array}
$$
It remains to verify the extendability of line bundles from $X$ to $V$. Fortunately this issue is addressed by the Barth-Lefschetz-type criterion~\cite[Theorem 2.2]{soms-vdv}, which implies that the restriction $\Pic(V)\to\Pic(X)$ is an isomorphism as soon as $\dta\les(\ell-2)/2$. Let us remark: 
$$
\dta\ges2\;\Rightarrow\;\ell\ges 6,\qquad\text{and}\qquad 
\bigg\lfloor\frac{\ell+1}{3}\bigg\rfloor\les\bigg\lfloor\frac{\ell-2}{2}\bigg\rfloor,\;\forall\;\ell\ges 6.
$$ 
Thus one should only impose that the differences $\eL_{j-1}^{-1}\eL_j$ are semi-ample; since $V$ is rational homogeneous, this is the same as being effective. 
 
\begin{m-theorem}\label{thm:homog}
\begin{enumerate}[leftmargin=4ex]
\item Let $V=G/P$ be a rational homogeneous variety, with $\ell\ges6$. Suppose $X$ is a smooth subvariety of codimension $\dta\les\lfloor(\ell+1)/3\rfloor$, whose normal bundle splits into a direct sum of line bundles which are `ordered': 
$$\eN_{X/V}= \ouset{j=1}{\dta}{\bigoplus}\eL_j\otimes \eO_X,
\quad\text{with}\; \eL_{j-1}^{-1}\eL_j\in\Pic(V)\;\text{effective, for}\; j=1,\dots,\dta.$$
Then $X$ is a complete intersection in $G/P$.  

\item 
If $G$ is simple and $P$ is maximal, the ordering condition is automatically satisfied, and the values of $\ell,\dta$ are as in the table below:
\renewcommand{\arraystretch}{1.5}\begin{m-eqn}{
\begin{tabular}{|c|c|c|c|c|c|c|}
\hline
{\rm group} $G$ & 
${\rm SL}(\ell+1)$ & 
${\rm SO}(2\ell+1)$& 
${\rm Sp}(\ell)$& 
${\rm SO}(2\ell)$& 
$E_6,\, E_7$ & $E_8$
\\ \hline 
\textbf{$\dta\les\ldots$} & 
$\big\lfloor\frac{\ell+1}{3}\big\rfloor,\,\ell\ges6$ &
$\big\lfloor\frac{\ell+1}{3}\big\rfloor,\,\ell\ges5$ &
$\big\lfloor\frac{\ell+1}{3}\big\rfloor,\,\ell\ges6$ &
$\big\lfloor\frac{\ell+1}{3}\big\rfloor,\,\ell\ges5$ &
$2$ & $3$
\\ \hline 
\end{tabular}
}\label{eq:ell1}\end{m-eqn}\renewcommand{\arraystretch}{1}
\end{enumerate}
\end{m-theorem}
In the second case, the reason for the changes in the table comes from the knowledge of the exact value of the partial amplitude of $\eT_V$ (cf.~\cite{gold}), which intervenes in the Barth-Lefschetz isomorphism (cf.~\cite[Theorem 2.2]{soms-vdv}, also~\cite{sena}). 

\begin{m-remark}
Typically, the embedding dimension of homogeneous varieties greatly exceeds their dimension. However, the same criterion for complete intersections as~\cite[Satz 5]{falt-krit} holds for subvarieties of codimension up to $\lfloor(\ell+1)/3\rfloor$. Similar statements hold for the isotropic (orthogonal and symplectic) Grassmannians (corresponding to the groups of type $B_\ell, C_\ell, D_\ell$). Note also that we obtain `exotic' criteria, \textit{e.g.} for $2$-codimensional subvarieties of the Cayley plane and Freudenthal's variety. 

On the lower side, we remark that for $\mbb P^\ell$ we get a weaker bound than $\ell/2$, as Faltings does. This is due to the fact that his proof uses a technical $G3$-criterion for subvarieties of projective spaces. In our treatment, this issue is hidden at the first step of the proof of Lemma~\ref{lm:surj}, since the partial amplitude assumption yields the $G3$-property (cf.~\cite[Section 2.2]{hlc}). 
\end{m-remark}


\subsection{On the number of equations defining subvarieties}\label{ssct:number}

In order to apply it, the previous theorem raises the question how to decide the splitting of the normal bundle (of the subvariety $X$ of the rational homogeneous variety $V$). The next result says that the normal bundle $\eN_{X/V}$ splits as soon as the ideal defining $X$ is generated by `few enough' elements. 
The case where the ambient space is $\mbb P^\ell$ is treated in~\cite[Satz 2, 3]{falt-krit}); the proof relies on a general lemma, recalled below. 

\medskip\nit\textbf{Hauptlemma}\quad \textit{%
Let $X$ be a projective variety and 
$$
0 \to\eE\srel{\alpha}{\to}\bbM\srel{\beta}{\to}\eF\to0,\quad\text{ with }\bbM=\ouset{j=1}{a}{\bigoplus}\eM_j,\quad\eM_j\in\Pic(X).
$$ 
be an exact sequence of vector bundles. (We say that ${\{\eM_j\}}_j$ generate $\eF$.) \newline
Take $s\in\Hom(\eM_j,\bbM)$ whose $j^{\rm th}$ component is the identity, and denote by $Z, W$ the zero loci of $\beta\circ s$ and ${\rm pr}_j\circ\alpha$, respectively (${\rm pr}_j$ stands for the projection onto $j^{\rm th}$ component). Then the intersection $Z\cap W$ is empty.%
}\medskip 

Now we specialize this result to $X\subset V=G/P$, of codimension $\dta$; the bounds below, on the number of generators, are in the same vein as in \lcit\ For shorthand, let 
$$
a:={\rm rank}(\bbM),\;f:={\rm rank}(\eF),\;e:={\rm rank}(\eE).
$$
We assume that the line bundles ${\{\eM_j\}}_j$ are ordered: 
\begin{m-eqn}{
\eM_{j-1}^{-1}\eM_j,\;j\ges2,\;\text{are effective.}\quad\text{(Equivalently, they are globally generated.)}
}\label{eq:order}\end{m-eqn} 

\begin{m-lemma}\label{lm:haupt}
{\rm(i)} Suppose that either one of the following inequalities hold: 
$$\;3\dta+f+a\les\ell+1\quad\text{or}\quad 3\dta+e+a\les\ell+1.$$ 
Then $\eF$ is isomorphic, through $\beta$, to the direct sum of $f$ line bundles in $\bbM$. (A similar statement holds for $\eE$.)

\nit{\rm(ii)} Suppose that $\eN_{X/V}^\vee$ (or $\eN_{X/V}$) is generated by $a\ges\dta$ line bundles (defined over $V$) which are ordered, such that the `number of relations' between the generators satisfies:
$$
\text{either}\quad (a-\dta)\les\ell+1-5\dta\qquad\text{or}\qquad  (a-\dta)\les\frac{\ell+1}{2}-2\dta.
$$
Then the normal bundle splits and it is isomorphic to the direct sum of $\dta$ of these line bundles.
\end{m-lemma}

\begin{m-proof} {\rm(i)} Let us assume that the first inequality holds. We argue by induction on the pair $(e,f)$, ordered lexicographically; if $e=0$, there is nothing to prove. The ordering assumption~\eqref{eq:order} implies that, for $s\in\Hom(\eM_1,\bbM)$ general, $\beta\circ s$ has one of the following properties: 
\begin{itemize}[leftmargin=4ex]
\item 
$Z=\emptyset.$
\item[] 
Then $\eF':=\eF/(\beta\circ s)\eM_1$ is locally free and we have $0{\to}\eE{\to}\ouset{j=2}{a}{\bigoplus}\eM_j\srel{\beta'}{\to}\eF'{\to}0.$ The induction hypothesis implies ($f$ decreases) that $\eF'$ splits (with splitting determined by $f-1$ components of $\beta'$; we abusively write $\beta'^{-1}$ for its inverse). Then $\eF$ splits too, because $\beta\beta'^{-1}$ yields an inverse to $\eF\to\eF'$.\medskip 
\item
$Z\subset X$ is smooth, non-empty, and $Z\cap W=\emptyset.$ 
\item[] 
In this case, we are going to show that $W=\emptyset$; that is, $\eE\to\eM_1$ is surjective. Hence the kernel $\eE'$ fits into $0{\to}\eE'{\to}\ouset{j=2}{a}{\bigoplus}\eM_j{\to}\eF{\to}0$ ($e$ decreases), so $\eF$ splits. 
\item[]
So let us prove that $W$ is empty. Faltings' bound applied to $Z\subset V$ yields 
$$
\cd(X\sm Z)\les\cd(V\sm Z)\les\dim V-\ell+2\delta_Z-2,\quad\delta_Z=\codim_V(Z)=\dta+f.
$$ 
(During the inductive process, this is decreasing with $f$.) 

If $W\neq\emptyset$, then it is a projective subscheme of $V$ of dimension at least $\dim X-e$. (During the inductive process, when $e$ decreases, this quantity increases.) Now observe that the hypothesis yields $\dim W\ges\cd(V\sm Z)+1$, a contradiction.
\end{itemize}
If the second inequality of the proposition is satisfied, we dualize the sequence and repeat the previous argument. (Here one starts with a general homomorphism $\eM_a^{-1}\to\bbM^\vee$.) 

\nit{\rm(ii)} In this case, we have $f=\dta$. 
\end{m-proof}

\begin{m-theorem}\label{thm:generators}
Let $ V=G/P$ be a rational homogeneous variety, with $\ell\ges 7$. Let $X$ be a $\dta$-codimen\-sional subvariety, $\dta\les\frac{\ell+1}{4}$, whose ideal is generated by `polynomials' $v_j\in\Gamma(V,\eM_j),$ $j=1,\dots,a$, such that the line bundles ${\{\eM_j\}}_j$ are ordered (cf.~\eqref{eq:order}). Then $X\subset V$ is a complete intersection, as soon as any of the following inequalities is satisfied: 
$$a\les \frac{\ell+1}{2}-\dta\quad\text{or}\quad a\les\ell+1-4\dta.$$
\end{m-theorem}
When $P$ is maximal parabolic, the Picard group of $V$ is cyclic and we obtain the statement in the Introduction.

\begin{m-proof}
Indeed, the differentials of the sections $\{v_j\}_j$ along $X$ yield a surjective homomorphism $\ouset{j=1}{a}{\bigoplus}\eM_j^{-1}\otimes\eO_X\to\eI_X/\eI_X^2=\eN_{X/V}^\vee$. It remains to apply Theorem~\ref{thm:homog} and Lemma~\ref{lm:haupt}.
\end{m-proof}


\section{On an issue raised by Hartshorne}\label{sct:hart}

In the remaining part of the article we shall focus the two-codimensional case. We are going to work over the ground field $\kk=\mbb C$. Keeping in mind Hartshorne's conjecture~\cite{hart-conj}, one would like to remove the assumption on the splitting of the normal bundle. In this sense, the references~\cite{bart+vdv,ran,bal+chi,holm,holm+sch} prove that subvarieties $X\subset\mbb P^n$ of sufficiently low degree are complete intersections; that is, the splitting of the normal bundle is automatically satisfied. Here we are going to show that several results still hold for homogeneous varieties.


\subsection{The framework}\label{ssct:setup}

Suppose $V=G/P$, where $G$ is simple of rank $\ell$ and $P$ is a maximal parabolic subgroup. The choices for $G$ and the corresponding values of $\ell$ are listed below: 
\renewcommand{\arraystretch}{1.5}\begin{m-eqn}{
\begin{tabular}{|c|c|c|c|c|c|}
\hline 
group $G$&${\rm SL}(\ell+1)$&${\rm SO}(2\ell+1)$&${\rm Sp}(\ell)$&${\rm SO(2\ell)}$&$E_6,\,E_7,\,E_8,\,F_4$
\\ 
\hline 
$\begin{array}{c}\Pic(V)\srel{\cong}{\to}\Pic(X)\\[-.5ex] \text{(See Lemma~\ref{lm:pic-iso}.)}\end{array}$
&$\ell\ges 6$&$\ell\ges4$&$\ell\ges6$&$\ell\ges5$&yes $\checkmark$
\\ \hline 
\end{tabular} 
}\label{eq:ell2}\end{m-eqn}\renewcommand{\arraystretch}{1} 
A basis of the singular (co)homology group of $V$ is given by the Schubert subvarieties; the incidences between them are encoded in the so-called Hasse-diagram. Note that $\Pic(V)=H^2(V;\mbb Z)=\mbb Z\cdot\eO_V(1)$ is determined by the simple root defining $P$. The generator $l$ of $H_2(V;\mbb Z)$ is represented by lines $L$ of degree one (with respect to $\eO_V(1)$), embedded in $V$; we call them \emph{straight lines}, in analogy with $\mbb P^n$. We denote 
$$
m=m(V):=c_1(V)\cdot l-3.
$$
(The notation is chosen to coincide with the number $m$ appearing in Ran's article.) For the next computations, we find convenient to have at hand the values of $m(V)$ and the positivity $p(V)$ of the tangent bundle of $V$. (See Definition~\ref{def:part-pos}.) We took the values of $c_1(V)$ from Snow~\cite[Corollary 2.4]{snow} and of $p(V)$ from Goldstein~\cite{gold}. 
\renewcommand{\arraystretch}{1.35}\begin{m-eqn}{
\scalebox{.63}{ \begin{tabular}{|c|c|c|c|c|c|c|c|c|c|c|c|}\hline 
$G$&${\rm SL}(\ell+1)$&\multicolumn{2}{c|}{${\rm SO}(2\ell+1)$}&\multicolumn{2}{c|}{${\rm Sp}(\ell)$}&\multicolumn{2}{c|}{${\rm SO(2\ell)}$}&$E_6$&$E_7$&$E_8$&$F_4$
\\ \hline 
\multirow{2}{*}{$V=\frac{G}{P}$}&$\Gr(k;\ell+1)$&{$\oGr(k;2\ell+1)$}&\multirow{2}{*}{$\oGr(\ell;2\ell+1)$}&\multicolumn{2}{c|}{$\sGr(k;2\ell)$}&$\oGr(k;2\ell)$&$\oGr(k;2\ell)$&$\ldots$&$\ldots$&$\ldots$&$\ldots$ 
\\ 
&$k\les\ell$&$k\les\ell-1$&&$2\les k \les\ell-1$&$k=\ell$&$k \les\ell-2$&$k=\ell-1,\ell$&&&&
\\ \hline 
$m(V)$&$\ell-2$& $2\ell-k-3$ & $2\ell-3$ &\multicolumn{2}{c|}{$2\ell-k-2$}&2$\ell-k-4$&$2\ell-5$&$\ldots$&$\ldots$&$\ldots$&$\ldots$
\\ \hline 
$p(V)$&$\ell$&{$2\ell-2$}&$2\ell-1$&\multicolumn{2}{c|}{$2\ell-k$}&\multicolumn{2}{c|}{$2\ell-3$}&$11$&$17$&$29$&$8,9,10$ 
\\ \hline 
{$\spp_V$}& $\text{max. of}$ & $\text{max. of}$ & \multirow{2}{*}{$\ell-1$} & $\text{max. of}$ & \multirow{2}{*}{$\ell-1$} & $\text{max. of}$ & \multirow{2}{*}{$\ell-1$} & \multirow{2}{*}{$\ges3$}& \multirow{2}{*}{$\ges4$}& \multirow{2}{*}{$\ges4$}&
\\ 
Def.~\ref{def:msp}. & $k,\ell-k+1$ & $k,\ell-k$ & & $k,\ell-k$& & $k,\ell-k$ & &&&&
\\ \hline 
\end{tabular} }
}\label{eq:mp}\end{m-eqn}\renewcommand{\arraystretch}{1}


Now consider a smooth, $2$-codimensional subvariety $X$ of $V$. Since $H^4(V;\mbb Z)$ is not necessarily cyclic, we assume that the Poincar\'e-dual class to $X$ is a multiple of $[\eO_V(1)]^2$: 
$$
[X]=d_X\cdot\chi^2,\quad\text{with}\;\chi:=[\eO_V(1)]\in H^2(V;\mbb Z). 
$$
This condition is certainly satisfied by complete intersections. We refer to $d_X$ as \emph{the degree} of $X$.

\begin{m-lemma}\label{lm:pic-iso}
For $G,\ell$ as in~\eqref{eq:ell2}, the restriction $\Pic(V)\to\Pic(X)$ is an isomorphism. 
\end{m-lemma}

\begin{m-proof}
The Barth-Lefschetz-type theorem (cf.~\cite[Theorem 2.2]{bart+vdv}, \cite{sena}) implies that the relative homotopy groups $\pi_j(V,X)$ vanish for $j\les 3\les p(V)-3$.
\end{m-proof}

The Hartshorne-Serre construction~\cite{arro} shows that $X$ is the vanishing locus of a rank-two vector bundle on $V$. So we have a locally free resolution  
\begin{m-eqn}{
0\to\eO_V(-n)\srel{\alpha}{\to}\eN^\vee\srel{}{\to}\eI_X\to 0,
}\label{eq:ONX}\end{m-eqn}
where $\eN$ is a rank-two vector bundle on $V$, with the following Chern classes: 
\begin{m-eqn}{
\begin{array}{r}
c_2(\eN)=[X]=d_X\cdot\chi^2,\quad\det\eN=\eO_V(n),\\[1.5ex] 
\eN\otimes\eO_X\cong\eN_{X/V},\quad c_1(\eN_{X/V})=\eO_X(n).
\end{array}
}\label{eq:nNX}\end{m-eqn}
The tangent sequence to $X$ yields $n=n_X=c_1(V)-c_1(X)$. Moreover, we remark that $n>0$, since $\eN_{X/V}$ is globally generated and partially ample (cf.~\eqref{eq:faltings}), being a quotient of $\eT_{V}$. 
For $k\in\mbb Z$, let us denote:  
$$
\begin{array}{rl}
e(k):=d_X-nk+k^2,&\quad\text{so}\;\;c_2(\eN^\vee(k))=e(k)\cdot\chi^{2}.
\\[1ex]
\Delta(\eN):=4d_X-n^2,&\quad\text{so}\;c_2\big(\cal E{\kern-1pt}nd(\eN)\big)=\Delta(\eN)\cdot\chi^{2},\quad\text{the discriminant.}  
\end{array}
$$
Below is represented the $(d_X,n_X)$-plane with the curves which are involved in our discussion.

\begin{m-eqn}{
\begin{tikzpicture}[domain=-1:10]
    \draw[->] (-.5,0) -- (9,0) node[right] {$d$};
    \draw[->] (0,-.5) -- (0,5) node[above] {$n$};
    
    \draw[name path=U, color=white, scale=.25, domain=0:9.5, variable=\d] plot ({\d}, {18});
    \draw[name path=UU, color=white, scale=.25, domain=9.5:20, variable=\d] plot ({\d}, {18});
    \draw[name path=UUU, color=white, scale=.25, domain=20:23, variable=\d] plot ({\d}, {18});
    \draw[name path=R, color=white, scale=.25, domain=0:10, variable=\n] plot ({28}, {\n});

    \draw[name path = D, thick, scale=.25, domain=0:9.5, smooth, variable=\d] plot ({\d}, {2*sqrt \d});
    \draw[name path = DD, thick, scale=.25, domain=9.5:20, smooth, variable=\d] plot ({\d}, {2*sqrt \d});
    \draw[name path = DDD, thick, scale=.25, domain=20:27, smooth, variable=\d] plot ({\d}, {2*sqrt \d});
    \draw[thick, scale=.25, domain=27:33, smooth, variable=\d] plot ({\d}, {2*sqrt \d}) 
    node[below] {$\qquad\boldsymbol{n = 2\sqrt{d},\;\Delta=0}$};

    \draw[name path = d, thick, scale=.25, domain=0:9.5, smooth, variable=\d] plot ({\d}, {4*sqrt \d/3});
    \draw[name path = dd, thick, scale=.25, domain=9.5:20, smooth, variable=\d] plot ({\d}, {4*sqrt \d/3});
    \draw[name path = ddd, thick, scale=.25, domain=20:27, smooth, variable=\d] plot ({\d}, {4*sqrt \d/3});
    \draw[thick, scale=.25, domain=27:33, smooth, variable=\d] plot ({\d}, {4*sqrt \d/3}) 
    node[below] {$\qquad\boldsymbol{n = 2\cos(t_\mx)\sqrt{d}}$};
        
    \draw[name path = E, thick, scale=.25, domain=0:9.5, smooth, variable=\d] plot ({\d}, {\d/3+3}) node[left] at (0,3) {$m$};
    \draw[name path = EE, thick, scale=.25, domain=9.5:20, smooth, variable=\d] plot ({\d}, {\d/3+3});    
    \draw[name path = EEE, thick, scale=.25, domain=20:27, smooth, variable=\d] plot ({\d}, {\d/3+3});    
    \draw[thick, scale=.25, domain=27:31, smooth, variable=\d] plot ({\d}, {\d/3+3})  
    node[above] {$\begin{array}{l}\boldsymbol{e(m)=0}\\ \boldsymbol{n = \frac{d}{m}+m} \end{array}$};
        
    \draw[scale=.25, very thick, domain=0:18, smooth, variable=\n] plot ({9.5}, {\n}) node[below] at (9.5,0) {$m^2$};
    
    \draw[name path = H, very thick, scale=.25, domain=0:9.5, smooth, variable=\d] plot ({\d},6.15);

    \draw[name path = RR, very thick, dashed, scale=.25, domain=0:6.5, smooth, variable=\n] plot (20,{\n});
    \draw[name path = DR, very thick, scale=.25, domain=6.4:9, smooth, variable=\n] plot (20,{\n}) node[below] at (24,0) {$d_{Thm.\,\ref{thm:d>}}=O(\ell^4)$};
    
    \tikzfillbetween[of=U and H, on layer=ft] {pattern=horizontal lines, thin};
    \tikzfillbetween[of=E and H, on layer=ft] {pattern=grid, thin};
    \tikzfillbetween[of=E and D, on layer=ft] {pattern=vertical lines, thin};
    \tikzfillbetween[of=EE and UU, on layer=ft] {pattern=horizontal lines, thin};
    \tikzfillbetween[of=EEE and UUU, on layer=ft] {pattern=horizontal lines, thin};
    \tikzfillbetween[of=DD and EE, on layer=ft] {pattern=checkerboard};
    \tikzfillbetween[of=EEE and ddd, on layer=ft] {pattern=checkerboard};
        
    \node[below left] at (0,0) {$0$};\node[left]at(0,1.6) {$2m$};
    \node at (3.55, 1.3) {$\boldsymbol{\Delta\;>\kern-1.55ex>\;0}$}; 
    \node at (5,3.5) {$\boldsymbol{e(m)<0}$};
    \node at (2.37,1.52) {$\bullet$}; \node at (2.35,2) {$\boldsymbol{(m^2,2m)}$};
    \node at (7.4,2.85) {$^{unknown}$};  
    \node at (7.3,1.95) {$^{region}$};   
    \node at (1.5, 0.35) {$empty$};    \node at (3.5, 0.35) {$region$}; 
\end{tikzpicture}
}\label{eq:dn-plane}\end{m-eqn}
The figure makes easier to summarise our results:
\begin{enumerate}[leftmargin=4ex]
\item 
For the $\Delta>0$ region, below the middle parabola, we prove that the possible points $(d_X,n_X)$ are actually confined in the `comet shaped region' between the two parabolas; thus there are no $(d,n)$-pairs realised by subvarieties $X$ below the lower parabola. Moreover, $d_X$ is bounded below by a constant, depending on $V$, so the point $(d_X,n_X)$ belongs to the checker-like region.

\item 
Concerning the $\Delta\ges0$ region, we show that subvarieties with $(d_X,n_X)$ in either the \begin{tikzpicture}\draw [pattern=horizontal lines] (0,0) rectangle (0.5,0.25);\end{tikzpicture} or \begin{tikzpicture}\draw [pattern=grid] (0,0) rectangle (0.5,0.25);\end{tikzpicture} or \begin{tikzpicture}\draw [pattern=vertical lines] (0,0) rectangle (0.5,0.25);\end{tikzpicture}-shaded regions have split normal bundle, so they are complete intersections. 

\item 
From these facts, we deduce that whenever the degree $d_X$ is at most $m(V)^2$, the subvariety $X\subset V$ is a complete intersection. 
\end{enumerate}


\subsection{$\boldsymbol{\Delta(\eN)>0}$, lower bounds for $\boldsymbol{d_X}$} 

The typical situation for this case is when the normal bundle is stable. 
We settle it by adapting Schneider's methods~\cite{sch} to our setting. The assumption implies that we can write: 
$$
n=\zeta+\bar\zeta=2r\cos(t),\quad d_X=\zeta\bar\zeta=r^2,\quad \zeta=r\cdot\exp(it),\;t\in(-\pi,\pi).
$$
In other words, the Chern roots of $\eN$ are imaginary, conjugate to each other. 

\begin{m-lemma}\label{lm:t}
We have the following estimate: 
$$|t|<\frac{\pi}{p(V)-1}=t_\mx.$$
Thus the point $(d_X,n_X)$ belongs to the `comet-shaped' region between two parabolas: 
$$1>\frac{n^2}{4d}>\cos(t_\mx)>1-\frac{5}{(p(V)-1)^2}\ges1-\frac{5}{(\ell-1)^2}.$$
\end{m-lemma}

\begin{m-proof}
By definition, the total Segre class is 
$$s(\eN_{X/V})=\frac{1}{(1-\zeta)(1-\bar\zeta)}=(1+\zeta+\zeta^2+\dots)(1+\bar\zeta+\bar\zeta^2+\dots),$$ 
which yields 
$$
s_j(\eN_{X/V})=\sum_{a=0}^{j}\zeta^a\bar\zeta^{j-a}=\frac{\zeta^{j+1}-\bar\zeta^{j+1}}{\zeta-\bar\zeta}=\frac{\sin(j+1)t}{\sin (t)}.
$$
The normal bundle of $X$ is a quotient of the tangent bundle of $V$, so it is still Sommese $q(V)$-ample. It follows that $p(\eN_{X/V})\ges\dim X-q(V)=p(V)-2$. 

Proposition~\ref{prop:segre} implies that, for $j\les p(\eN_{X/V})$, the Segre class $s_j(\eN_{X/V})$ is a strictly positive multiple of $\chi^j$, so $\sin(j+1)t>0$ for $j=0,\dots,p(V-2)$. It follows that $|t|<t_\mx$, so $\cos(t)>\cos(t_\mx)$.
\end{m-proof}

The Lemma immediately yields a universal lower bound for $d_X$, thus further restricting the values it may take. 

\begin{m-definition}\label{def:msp}
Let $\spp_V$ be the dimension of  a maximal linear subspace of $V$ (Schubert subvariety, isomorphic to a projective space, and maximal with this property). We denote by $\msp_V$ such a `maximal Schubert projective subspace of $V$'; note that it is not necessarily unique. 

The values of $\spp_V$ for the classical homogeneous spaces are tabulated in the bottom row of the Table~\ref{eq:mp}. 
\end{m-definition}

\begin{m-theorem}\label{thm:d>}
Let $X\subset V$ be such that $\Delta(\eN)>0$. Then the degree of $X$ satisfies the inequality:
$$ 
d_X>\frac{\spp_V^2}{240}\cdot(p(V)-1)^2.
$$
That is, $(d_X,n_X)$ is in the \begin{tikzpicture}\draw [pattern=checkerboard] (0,0) rectangle (0.5,0.22);\end{tikzpicture}-shaded region. 
\\ Consequently, if $d_X$ does not satisfy this inequality, then necessarily $\Delta(\eN)\les0$. 
\end{m-theorem}

\begin{m-proof} 
Indeed, we have: 
$$
4d_X\cdot\frac{\pi^2}{(p(V)-1)^2}\ges 4r^2\sin^2(t)=\Delta(\eN)=\Delta(\eN\otimes\eO_{\msp_V})\ges \frac{\spp_V^2}{6},
$$
where, for the last estimate, we used Proposition~\ref{prop:sch}. 
\end{m-proof}

\begin{m-remark}\label{rmk:delta}
\begin{enumerate}[leftmargin=4ex]
\item 
We observe that the lower bound for $d_X$ is of \emph{order four} in $\ell$. This is especially interesting for the classical groups of type $A_\ell-D_\ell$. 
\item 
Let us consider $V=\mbb P^\ell$. The best previously known bound~\cite[Remark 4.7]{holm+sch} is $(\ell-1)(\ell+5)$, while our Theorem yields 
$$
d_{X\subset\mbb P^\ell}>\frac{\ell^2(\ell-1)^2}{24\pi^2}>\frac{\ell^2(\ell-1)^2}{240}.
$$ 
For $\ell\ges18$, our estimate is better. 
\item 
For low values of $\ell$, one can improve the lower bound by using the explicit minimal values for $\Delta(\eN\otimes\eO_{\msp_V})$ (cf.~\cite[Table p.~155]{sch}) rather than Proposition~\ref{prop:sch}. We shall use this in the proof of Theorem~\ref{thm:hart}.
\item 
But the formula is not restricted to projective spaces. For instance, in the case of the Grassmannian $V=\Gr(k;\ell+1),\;k\ges(\ell+1)/2$, one obtains $\frac{k^2(\ell-1)^2}{240}.$ 
\end{enumerate}
\end{m-remark}


\subsection{$\boldsymbol{\Delta(\eN)\les0}$, on a result of Ran}\label{ssct:ran}

Let the situation be as in Section~\ref{ssct:setup}. 
We are going to show that the argument of Ran~\cite{ran} can be generalized. The only change is that we find more comfortable working with stable maps (in the sense of Kontsevich-Manin) instead of secants, in order to understand the various strata which appear. 

Let $\tld V\srel{\sigma}{\to} V$ be the blow-up along $X$; the exact sequence~\eqref{eq:ONX} implies that $\tld V\subset\mbb P(\eN^\vee)=Proj(\Sym^\bullet\eN^\vee)$. The exceptional divisor is $E_X=\mbb P(\eN_{X/V}^\vee)\srel{\sigma_X}{\to} X$; it is a $\mbb P^1$-bundle over $X$. The fibres of $\sigma_X$ are called \emph{vertical lines} and their homology class is denoted $\phi$. Note that 
$$
H_2(\tld V;\mbb Z)= H_2(V;\mbb Z)\oplus\mbb Z\cdot\phi=\mbb Z\cdot l\oplus\mbb Z\cdot\phi.
$$ 
The lifting of $l\in H_2(V;\mbb Z)$ to $H_2(\tld V;\mbb Z)$ is given by the pre-image of a straight line which avoids $X$. 

\begin{m-definition}
\begin{enumerate}[leftmargin=4ex]
\item 
For $k\ges 0$, we denote $$b_k:=l-k\phi\in H_2(\tld V;\mbb Z).$$ It is determined by the conditions $\sigma_*(b_k)=l$ and $b_k\cdot E_X=k$. 

\item 
Let $\ovl M_{0,p}(\tld V;b_k)$ be the moduli space of $p$-pointed, genus zero stable maps, which represent the homology class $b_k$ (cf.~\textit{e.g.}~\cite{fult+pand}). The latter are denoted $(T,\ut, u)\srel{u}{\to}\tld V$, so $T$ is a quasi-stable curve (tree) whose components are isomorphic to $\mbb P^1$, $\ut=(t_0,t_1,\dots)$ stands for the markings, and $u:T\to\tld V$ is a morphism. 

\item[] 
The reason for introducing $b_k$ is clear: if a straight line in $V$ meets $X$ at $k$ points with multiplicity one, its proper transform in $\tld V$ represents $b_k$. Conversely, any $T\srel{u}{\to} \tld V$ representing $b_k$ has a unique `$\main$' component $T_{\main}\subset T$ isomorphic to $\mbb P^1$, which is mapped by $\sigma u$ to a straight line in $V$. The other `$\phi$-components' of $T$ are sent to vertical lines, thus $u_*[T_{\main}]=b_j,\,j\ges k.$ 

\item 
Let $\pnt\in V\sm X$ be a general point, in a sense yet to be made precise (cf.~Lemma~\ref{lm:empty}). For $k\ges 1$, we define: 
$$
\Sigma_k:=\ovl M_{0,1}(\tld V;b_k)\times_{\tld V}\{\pnt\}\subset\ovl M_{0,1}(\tld V;b_k).
$$ 
It is the analogue of the space of $k$-secants to $X$, passing through $\pnt$; by abuse of language, we call it the same. 
\item[] 
Note that $\ovl M_{0,1}(\tld V;\phi)\cong E_X\to X\cong\ovl M_{0,0}(\tld V;\phi)$, thus there are gluing morphisms 
$$
\big(\,\ovl M_{0,2}(\tld V;b_{k+1})\times_{\tld V}\{pt\}\,\big)\times_{\tld V} E_X
\;\srel{g_k}{\to}\;\ovl M_{0,1}(\tld V;b_k)\times_{\tld V}\{pt\}=\Sigma_k,
$$ 
which insert a vertical line (a `bubble') at the second marked point. 
\item[]
The term $\ovl M_{0,2}(\tld V;b_{k+1})\times_{\tld V}\{pt\}$ can be thought of as the universal curve over $\Sigma_{k+1}$. 
\item[] 
We say that the (schematic) image of $g_k$ is \emph{the boundary} of $\Sigma_k$, denoted $\partial\Sigma_k$.
\end{enumerate}
\end{m-definition}

It is known that $M_{0,0}(\tld V;b_0)\cong M_{0,0}(V;l)$ is irreducible, normal, of dimension $c_1(V)\cdot l+\dim(V)-3$, so we have 
$$
\Sigma_0\neq\emptyset,\quad\dim(\Sigma_0)=c_1(V)\cdot l-2=m+1.
$$
Since $X$ intersects any closed subscheme of $V$ of dimension at least two, we have $\Sigma_1\neq\emptyset$. As $k$ increases, at each step the dimension of the irreducible components of $\Sigma_k$ should drop by one, since the main component has to intersect $E_X$ at one more point or with higher multiplicity. 
Indeed its expected dimension (given by the relevant obstruction theory) is: 
$$
\mathrm{exp.dim.}(\Sigma_k)=c_1(\tld V)\cdot b_k+\dim V-2-\dim V=m+1-k.
$$ 
The last non-empty secant variety should be $\Sigma_{m+1}$, a finite set of points. However, this is not necessarily the case, and we are interested when does the process terminate earlier, after at most $m$ steps. In such a situation, at some moment the dimension drops by two, at least.

\begin{m-lemma}\label{lm:jump}
Suppose that $\codim_{\Sigma_k}\partial\Sigma_{k}\ges2.$ Then we have  $e(k)=0$. 
\\(The codimension is the maximal codimension of the irreducible components.) 
\end{m-lemma}

\begin{m-proof}
We observe that, for any $(T, t_0, u)\in\Sigma_k$, the marking $t_0$ is mapped to $pt\not\in E_X$, hence $t_0$ belongs to $T_\main$; it's the only component whose image is not contained in $E_X$. For this reason one can identify $(T,t_0,u)$ with $(T,u)$, by means of $\Sigma_k\to\ovl M_{0,0}(\tld V;b_k)$.
 
We claim $(T,t_1,u)$ belongs to $\partial\Sigma_k$ if and only if $T\neq T_\main$. The necessity is clear, let us prove the sufficiency. If $T$ contains other components beside $T_\main$, then there are non-constant vertical components, due to the stability of the map. It follows that $u_*[T_\main]=b_j,\;j\ges k+1$, so $(T,u)$ is obtained by inserting bubbles to $(T_\main,u)\in\Sigma_j$. 

The remaining part of the proof is similar to~\cite[Proposition, p.~334]{ran}. The hypothesis implies that there is a smooth, projective curve $C$ and a morphism 
$C\to \Sigma_k\sm \partial \Sigma_{k}.$ 
In other words, there is a complete, $1$-dimensional family of stable maps $(\mbb P^1,u)$ passing through $pt$ and representing $b_k$. Let $S$ be the corresponding universal curve; that is,   
$$
S:=C\times_{\ovl M_{0,1}(\tld V;b_k)}\ovl M_{0,2}(\tld V;b_k).
$$ 
It is a ruled surface over $C$ which fits in the diagram 
$$
\xymatrix@R=2em{
S\ar[d]^-\varpi\ar[r]^-u&\tld V
\\ 
C\ar@/^3ex/[u]^-{t_0}&
}
$$
Here $u$ is the evaluation morphism at the second marked point and $t_0$ is the section given by the (first) marked point which is mapped to $pt\in V\sm X=\tld V\sm E_X$.

Now recall that $\tld V\subset\mbb P(\eN^\vee)$. The restriction to $\tld V$ of the Euler sequence on $\mbb P(\eN^\vee)$ is
$$
0\to \eO_{\tld V}(E_X)\srel{\tld\alpha}{\to} \sigma^*\eN^\vee\to\tld\eQ\to0,\quad (\tld Q\;\text{is the universal, locally free quotient}).
$$
Note that $\tld\alpha$ is actually induced by the pull-back of $\alpha$ in~\eqref{eq:ONX}; indeed, the latter vanishes precisely along $E_X$. 
We further pull-back this sequence to $S$, and we wish to compute $c_1(u^*\eO_{\tld V}(E_X))=u^*c_1(\eO_{\tld V}(E_X))$. For $f:=[\text{fibre of}\;\varpi]\in H_2(S;\mbb Z)$, we have 
$$
c_1(u^*\eO_{\tld V}(E_X))\cdot f=E_X\cdot b_k=k=k\cdot (u^*\sigma^*\chi\cdot f),\quad (\chi=[\eO_V(1)]),
$$ 
so $c_1(u^*\eO_{\tld V}(E_X))= k\cdot u^*\sigma^*\chi+k'\cdot f,$ with $k'\in\mbb Z$. (This is because $S$ is a ruled surface.) We compute $k'$ by restricting to $t_0(C)$; it is contracted by $u$, so $k'$ vanishes. We conclude that: 
\\ \null\hfill$\begin{array}[b]{l}
c_1(u^*\eO_{\tld V}(E_X))=k\cdot u^*\sigma^*\chi,\quad c_1(u^*\tld\eQ)=(d-k)\cdot u^*\sigma^*\chi,  
\\ d_X\cdot u^*\sigma^*\chi^2=c_2(u^*\sigma^*\eN^\vee)=k(d-k)\cdot u^*\sigma^*\chi^2.
\end{array}$\hfill\null 
\end{m-proof}

\begin{m-lemma}\label{lm:empty}
If $X$ is contained a hypersurface $Y\in|\eO_V(k)|$, then we have $\Sigma_j=\emptyset, \forall j\ges k+1$. 
\\ (The point $\pnt\in V$ involved in defining $\Sigma_0$ is chosen in the complement of $Y$.) 
\end{m-lemma}

\begin{m-proof}
Indeed, let $\tld Y=\sigma^*Y-\epsilon E_X,\,\epsilon\ges1,$ be the proper transform of $Y$ and take $u:T\to\tld V$ representing the class $b_j$. Since $u(T_\main)$ is not contained in $Y$ and represents $b_{j'},\;j'\ges j$, we have 
$\;0\les u_*[T_\main]\cdot\tld Y=(l-j'\phi)\cdot(k\chi-\epsilon E_X)=k-\epsilon j'.$ 
\end{m-proof}

\begin{m-lemma}\label{lm:D}
Suppose $\Delta(\eN)\les 0$. Then there is an integer $1\les k\les n/2$, such that $e(k)\ges0.$ Moreover, $X$ is contained in a hypersurface of degree $k$ and therefore we have $\Sigma_{k+1}=\emptyset$.
\end{m-lemma}

\begin{m-proof}
We distinguish two possibilities.\smallskip

\nit\underbar{$\eN$ is semi-stable}\quad  In this case, a theorem of Bogomolov implies that $\Delta(\eN)\ges0$. It follows that $\Delta(\eN)=0$, so we have $n=2k,\,c_1(\eN^\vee(k))=0,\,c_2(\eN^\vee(k))=0$. Thus $\eN(-k)$ has vanishing Chern classes and it is semi-stable and admits a non-trivial section. Otherwise $\eN(-k)$ is stable, so it possesses a flat Hermite-Einstein metric (by the Kobayashi-Hitchin correspondence), hence it is determined by a representation of the fundamental group of $V$, which is trivial; a contradiction, since we assumed that $\eN(-k)$ had no section. 

So there is a non-trivial homomorphism $\eO_V\to\eN(-k)$; its cokernel is locally free in codimension one. (Otherwise there was a non-trivial section of $\eN(-k-1)$, contradicting the semi-stability property.) Since $c_2(\eN(-k))=0$, the second Chern class of the cokernel vanishes, so it is actually an invertible sheaf with vanishing first Chern class. Therefore we have the exact sequence 
$$
0\to\eO_V\to\eN(-k)\to\eO_V\to0\quad\Rightarrow\quad \eN=\eO_V(k)^{\oplus 2}.
$$ 

\nit\underbar{$\eN$ is not semi-stable}\quad Let $k_\mn$ be the smallest integer $k$ such that there is a non-trivial homomorphism $\eO_V(-k)\to\eN^\vee;$ that is, a section $\nu$ of $\eN^\vee(k_\mn)$. It has the following properties: 
\begin{itemize}[leftmargin=4ex]
\item $k_\mn< n/2$\quad It is the non-semi-stability condition. In fact, $\eO_V(-k_\mn)$ is the first term of the Harder-Narasimhan filtration of $\eN^\vee$; this uniquely characterises $k_\mn$.  
\item $e(k_\mn)\ges0$\quad The minimality of $k_\mn$ implies that the vanishing locus of $\nu$ is a two-codimensional, possibly empty, subscheme which represents $c_2(\eN^\vee(k_\mn))$. 
\item $k_\mn\ges1$\quad Since $n\ges1$, the left-hand side of the exact sequence below vanishes
$$
0\to\Gamma(\eO_V(k_\mn-n))\to\Gamma(\eN^\vee(k_\mn))\to\Gamma(\eI_X(k_\mn))\to0,
$$
hence we obtain a non-vanishing section $\eI_X(k_\mn)\subset\eO_V(k_\mn)$.
\end{itemize}\medskip 
In both cases, the variety $X$ is contained in the hypersurface $Y$ of degree $k_\mn$, corresponding to the induced section in $\eO_V(k_\mn)$, so $\Sigma_j=\emptyset$ for $j>k_\mn$. 
\end{m-proof}

\begin{m-theorem}\label{thm:ran}
Let the situation as in~\ref{ssct:setup}, with $\Delta(\eN)=d-4n^2\les0$. 
\\ Assume that either one of the following conditions is satisfied (see Figure~\eqref{eq:dn-plane}): 
\begin{enumerate}[leftmargin=4ex]
\item $\frac{d_X}{m}+m\les n$\quad (the \begin{tikzpicture}\draw [pattern=horizontal lines] (0,0) rectangle (0.5,0.20);\end{tikzpicture} and \begin{tikzpicture}\draw [pattern=grid] (0,0) rectangle (0.5,0.20);\end{tikzpicture}-shaded regions);
\item $n/2\les m$\quad (the \begin{tikzpicture}\draw [pattern=vertical lines] (0,0) rectangle (0.5,0.21);\end{tikzpicture} and \begin{tikzpicture}\draw [pattern=grid] (0,0) rectangle (0.5,0.22);\end{tikzpicture}-shaded regions).
\end{enumerate}
Then $X\subset V$ is a complete intersection.  
\end{m-theorem}

\begin{m-proof} 
The hypotheses are made in such a way to ensure that $k_\mn\les m$: for (i), this follows from $e(k_\mn)\ges0\ges e(m)$, combined with the monotonicity of the function $j\mapsto e(j)$; for (ii), this is Lemma~\ref{lm:D}. Since $\Sigma_{k_\mn+1}=\emptyset$, the discussion preceding Lemma~\ref{lm:jump} implies that there is an $j\les m$ such that $\codim_{\Sigma_{j}}(\Sigma_{j+1})\ges2$; let $j_\mn$ be minimal with this property, so $1\les j_\mn\les k_\mn\les n/2.$ The monotonicity of the function $j\mapsto e(j)$ yields 
$$
0=e(j_\mn)\ges e(k_\mn)\ges0\quad\Rightarrow\quad j_\mn=k_\mn,\;\; e(k_\mn)=0.
$$
We conclude that $c_2(\eN^\vee(k_\mn))=0$, so the zero locus of the section $\nu$ is empty. One gets an exact sequence of vector bundles 
$$
0\to\eO_V(-k_\mn)\to\eN^\vee\to\eO_V(k_\mn-n)\to0,
$$
which splits, since the relevant ${\rm Ext}^1$-group vanishes. Thus $\eN$ splits, so $\eN_{X/V}$ too. By the second part of Lemma~\ref{lm:XV} (alternatively, Theorem~\ref{thm:g3}), $X$ is a complete intersection.  
\end{m-proof}


\subsection{Application to Hartshorne's problem}\label{ssct:hart}

Let the notation be as before. 

\begin{m-theorem}\label{thm:hart}
The subvariety $X\subset V$ is a complete intersection in any of the following cases: 
\begin{enumerate}[leftmargin=4ex]
\item $\ell\ges11$ and $d_X\les m(V)^2$;
\item $6\les\ell\les10$, $\,d_X\les3m(V)^2/10$, and also $V\neq \oGr(3;12), \oGr(3;13),\sGr(3;12)$.
\end{enumerate} 
\end{m-theorem}

\begin{m-proof}
The statement is clear from the picture~\eqref{eq:dn-plane}; in fact, we must prove that the dashed vertical line is indeed at the right of $\{d=m^2\}$. We distinguish several cases: 
\begin{itemize}[leftmargin=4ex]
\item 
$\Delta\les0$:\quad 
Here we may have $2m\les n$ or $n\les 2m$. The conclusion follows from the first, respectively the second part of Theorem~\ref{thm:ran}. 
\item 
$\Delta>0$:\quad 
Theorem~\ref{thm:d>} implies that this situation does not occur for $d_X\les\frac{\spp_V^2(p(V)^2-1)}{240}$. The table~\eqref{eq:mp} shows that $p(V)-1\ges m(V)$ and $\spp_V\ges\ell/2$. Thus, for $\ell\ges31$, we have $\frac{\spp_V^2(p(V)^2-1)}{240}\ges m(V)^2$. 

For the range $11\les\ell\les30$, we use Remark~\ref{rmk:delta}(iii). The proof of Theorem~\ref{thm:d>} shows that the degree of $X$ satisfies the inequality: 
$$
d_X>\frac{\Delta(\eN\otimes\eO_{\msp_V})}{4\pi^2}\cdot(p(V)-1)^2\ges \frac{\Delta_\mn(\spp_V)}{4\pi^2}\cdot m(V)^2,
$$
where $\Delta_\mn(\spp_V)$ denotes the minimal value of the discriminant on a $\spp_V$-dimensional projective space (cf.~\cite[Table, p.~155]{sch}). Since $\ell\ges11$ and $\spp_V\ges\lceil\ell/2\rceil$, it follows that $\spp_V\ges 6$ and $\Delta_\mn(\spp_V)\ges71$. In remains to use this estimate in the inequality above. 
\end{itemize}
For the second statement, the case $\Delta\les0$ remains unchanged. For $\Delta>0$ and $\ell=\{6,\dots,10\}$, we observe that $\spp_V\ges 4$ (with the three exceptions) and use that $\Delta_\mn(4)=12$.
\end{m-proof}


\appendix

\section{Partial positivity}\label{sct:partialpos}

For the reader's comfort, we recall here a few facts about partial ampleness of vector bundles and subvarieties. 

\begin{m-definition}\label{def:part-pos} 
Let $W$ be a projective variety. 
\begin{enumerate}[leftmargin=4ex]
\item (cf.~\cite{dps}) 
A \emph{flag} (`Babylonian tower') of normal and irreducible subvarieties of $W$ is a chain $W=Y_0\supset Y_1\supset\dots\supset Y_d,$ such that each $Y_j\subset Y_{j-1}$ is a Cartier divisor. Then $d$ is the \emph{length} of the flag and $Y_d$ is its \emph{end}.\smallskip  
\item (cf.~\cite{tota}) 
For $\eL\in\Pic(W)$ and an integer $0\les q\les\dim W$, we say that $\eL$ is \emph{$q$-ample} if, for any coherent sheaf $\eG$ on $W$, the following cohomology vanishing holds:
$$
\exists m_\eG\in\mbb Z\;\forall m\ges m_\eG\;\forall t>q,\quad H^t(W,\eG\otimes\eL^m)=0.
$$
(For $q=0$, $\eL$ is ample in the usual sense, while for $q=\dim W$ the condition is vacuous.)

\item 
A line bundle which is both $q$-ample and semi-ample ---some tensor power is globally generated--- is called \emph{$q$-positive}.
\item (cf.~\cite{arap}) 
A locally free sheaf $\eE$ on $W$ is \emph{$q$-ample} if $\eO_{\mbb P(\eE)}(1)$ on $\mbb P(\eE):=Proj(\Sym^\bullet_{\eO_Y}\eE)$ is $q$-ample. It is equivalent saying that, for any coherent sheaf $\eG$ on $W$, there is $m_\eG>0$ such that: 
$\,H^t(W,\eG\otimes\Sym^a(\eE))=0,\;\forall t>q,\; \forall a\ges m_\eG.$ 
\item[] 
The \emph{$q$-amplitude of $\eE$}, denoted by $q(\eE)$, is the smallest integer $q$ with this property. 
\item[] 
We call $p(\eE):=\dim W-q(\eE)$ the \emph{positivity} of $\eE$. (Goldstein~\cite{gold} calls the positivity \emph{co-ampleness}.)
\end{enumerate}
\end{m-definition}

\begin{m-proposition}{\rm (cf.~\cite{soms})} \label{prop:soms} 
For a \emph{globally generated}, locally free sheaf $\eE$ on $W$, the following statements are equivalent: 
\begin{enumerate}[leftmargin=5ex]
\item[\rm(i)] 
$\eE$ is $q$-ample; 
\item[\rm(ii)] 
In the diagram below, the fibres of the morphsm $\phi$ are at most $q$-dimensional:
$$\xymatrix@R=2em{
\mbb P(\eE)\ar[r]^-\phi\ar[d]_-\psi & \mbb P:=|\eO_{\mbb P(\eE)}(1)|
\\ 
W&
}$$
\end{enumerate}
We say that $\eE$ is \emph{Sommese-$q$-ample} if it satisfies any of these conditions.
\end{m-proposition}

Let us recall that the total Segre class of a vector bundle $\eF$ is by defined as $1/c(\eF^\vee)$, the inverse of the total Chern class of $\eF^\vee$, and the Segre classes are its graded components. 

\begin{m-proposition}\label{prop:segre}
Let $\eE$ be a globally generated vector bundle on a projective variety $W$, which is Sommese $q(\eE)$-ample. Then the first $p(\eE):=\dim X-q(\eE)$ Segre classes of $\eE$ are non-zero. 
\end{m-proposition}
The statement generalizes the well know fact~\cite[Theorem 1.5]{gies} saying that the Segre classes of ample vector bundles are strictly positive. 

\begin{m-proof}
We use the notation in the diagram above. Note that, for any fibre $F$ of $\phi$, the morphism $\psi_F:F\to W$ is finite, because the restriction of $\phi$ to the fibres of $\psi$ are closed embeddings. The Segre classes of $\eE$ have the following property: for any $j$-dimensional subvariety $Z\subset W$, one has ($e:=\mathrm{rank}(\eE)$):
$$
s_j(\eE)\cdot [Z]=[\eO_{\mbb P(\eE)}(1)]^{j+e-1}\cdot[\psi^*Z]=[\eO_{\mbb P}(1)]^{j+e-1}\cdot\phi_*[\psi^*Z].
$$
By hypothesis, the general fibre $F$ of $\phi$ is at most $q(\eE)$-dimensional. Let $Z\subset W$ be a general complete intersection of dimension $j\les p(\eE)$; it has the property that the intersection with $\psi(F)$ is either empty or finite. Thus the subscheme $\phi_*(\psi^*Z)\subset\mbb P$ is $(j+r-1)$-dimensional, and consequently the intersection product $s_j(\eE)\cdot Z$ is strictly positive. 
\end{m-proof}

\begin{m-definition} (cf.~\cite{hlc})\quad  
Let $Y\subset W$ be a closed subscheme, $\tld W$ the blow-up of the sheaf of ideals $\eI_Y\subset\eO_W$, and $E_Y$ the exceptional divisor. One says that $Y$ is \emph{$q_Y$-ample}, for $0\les q_Y\les\dim Y$, if $\eO_{\tld W}(E_Y)$ is $\tld q_Y=(q_Y+\codim_W(Y)-1)$-ample: for any coherent sheaf $\tld\eG$ on $\tld W$ there is an integer $m_{\tld\eG}$, such that it holds:
$$
\forall m\ges m_{\tld\eG},\;\forall t\ges\tld q_Y,\qquad H^t(\tld W,\tld\eG\otimes\eO_{\tld V}(mE_Y))=0.
$$
\end{m-definition}

\begin{m-proposition}{(cf.~\cite[Proposition~1.4]{hlc})}\label{prop:XV}
Suppose $W$ is smooth and $Y$ is a local complete intersection (lci, for short). Then $Y$ is $q_Y$-ample if and only if 
\begin{m-eqn}
{\eN_{Y/W} \;\text{is $q_Y$-ample}\quad\text{and}\quad \cd(W\sm Y)\les q_Y+\codim_W(Y)-1.}\label{eq:N}
\end{m-eqn}
Here $\cd(\cdot)$ stands for `cohomological dimension'.
\end{m-proposition}


\section{Schwartzenberger's conditions}\label{sct:schw}

Schneider used the integrality conditions (cf. Schwartzenberger~\cite[Theorem 22.4.1]{hirz}) imposed by the Riemann-Roch theorem in order to deduce lower bounds (cf.~\cite[Table pp. 155]{sch}) for the discriminant of rank two vector bundles over projective spaces. 
He obtained the table by using a computer and, based on this numerical evidence, claimed that the lower bound for the discriminant should increase quadratically with the dimension of the projective space; however, the proof ``leads to difficult number-theoretic questions''. In what follows, we show that such a quadratic lower bound holds indeed. 

Let $\eF$ be a rank-$2$ vector bundle on $\mbb P^p$ and $\dta_1,\dta_2$ its Chern roots. By the Riemann-Roch theorem, the Euler characteristic of $\eF(k), k\in\mbb Z,$ equals: 
$$
RR(p)=\binom{\dta_1+k+p}{p}+\binom{\dta_2+k+p}{p},
$$
so this quantity, which \textit{a priori} is a fraction with denominator $p!$, is actually an integer. To approximate factorials, we are going to use the classical Stirling formula (inequality):
$$
p!\ges \sqrt{2\pi p}\cdot\left( \frac{p}{\textbf e} \right)^p,
$$ 
where $\textbf e =2.718\dots$ is the basis of natural logarithms. 

We are interested in $\Delta(\eF)$, which is invariant under the twisting of $\eF$ by a line bundle, so we may assume that $c_1(\eF)$ equals $0$ or $1$; let 
$$
\dta:=\dta_1=-\dta_2,\quad\text{respectively,}\quad\dta_1=\dta+1/2,\;\dta_2=-\dta+1/2.
$$

\begin{m-proposition}\label{prop:sch}
Suppose $p\ges4$ and $\Delta(\eF)>0$. Then we have the following lower bound: 
$$
\Delta(\eF)> \frac{p^2}{6}.
$$
\end{m-proposition}

\begin{m-proof}
We start with the case when $p$ is even, so we replace $p$ by $2p$ and put $k=-p$. As we are going to see, this shift simplifies the expression of the Riemann-Roch number; it is the main difference compared to Schneider's approach, who considered $k=0$.

\smallskip\nit\unbar{Case $c_1(\eF)=0$}\quad So we are interested in $RR(2p):=\binom{\dta + p}{2p}+\binom{-\dta+p}{2p}$. We compute 
$$
\begin{array}{rl}
\binom{\dta + p}{2p}
&
=\frac{(\dta+p)(\dta+(p-1))\dots(\dta+1)\dta(\dta-1)\dots(\dta-(p-1))}{(2p)!}
\\[2ex] 
&
=\frac{(\dta^2+p\dta)(\dta^2-(p-1)^2)\dots(\dta^2-1)}{(2p)!}
=(-1)^{p}(d-p\dta)\cdot\frac{(d+(p-1)^2)\dots(d+1)}{(2p)!},
\end{array}
$$ 
and deduce that 
$$(-1)^{p}RR(2p)=2\cdot\frac{(d+(p-1)^2)\dots(d+1)d}{(2p)!}$$ 
is a strictly positive integer. By the classical inequality of means, we have 
$$
(d+(p-1)^2)\dots(d+1)d\les\bigg(d+\frac{(p-1)(p-1/2)}{3}\bigg)^p<\bigg(d+\frac{p^2}{3}\bigg)^p.
$$
By using Stirling's formula, we deduce:
$$
d+\frac{p^2}{3}>\sqrt[p]{\frac{(2p)!}{2}}>\frac{\sqrt[2p]{4\pi p}}{\sqrt[p]{2}}\cdot\frac{(2p)^2}{\mathbf e^2}>\frac{p^2}{2}
\quad\Rightarrow\quad d>\frac{p^2}{6}\quad\Rightarrow\quad\Delta>\frac{(2p)^2}{6}.
$$

\smallskip\nit\unbar{Case $c_1(\eF)=1$}\quad   Now the following quantity is an integer: 
$RR(2p):=\binom{\dta + p+1/2}{2p}+\binom{-\dta+p+1/2}{2p}$. 
$$
\begin{array}{rl}
\binom{\dta+p+1/2}{2p}
&
=\frac{(\dta + p+1/2)\dots(\dta+1/2)(\dta-1/2)\dots(\dta-p+3/2)}{(2p)!}
\\[2ex] 
&
=\frac{(\dta^2+2p\dta+ p^2-1/4)(\dta^2-(p-3/2)^2)\dots(\dta^2-(1/2)^2)}{(2p)!}.
\end{array}
$$ 
By replacing $\dta^2=-d+1/4$, we deduce that 
$$
(-1)^{p-1}RR(2p) =2\cdot\frac{(-d+p^2)\ouset{k=1}{p-1}{\prod}(d+(k-1/2)^2-1/4)}{(2p)!}
$$
is an integer. If $d\ges p^2>\frac{(2p)^2}{6}$, there is nothing to prove. 

So let us consider the case when $d<p^2$. The inequality of means implies that the Riemann-Roch number is bounded above by 
$\Big(\frac{p-2}{p}d+\frac{p^2+2}{3}\Big)^p.$ By using the Stirling formula, this yields the inequality 
$$
\begin{array}{r}
\disp 
\frac{p-2}{p}d+\frac{p^2+2}{3}\ges\frac{\sqrt[2p]{4\pi p}}{\sqrt[p]{2}}\cdot\frac{(2p)^2}{\mathbf e^2}>\Big(1+\frac{\ln(3p)}{2p}\Big)\frac{(2p)^2}{\mathbf e^2}
\\[2ex]
\disp\Rightarrow\quad 
d>\Big(\frac{4}{\textbf e^2}-\frac{1}{3}\Big)p^2+\frac{p\ln(3p)}{4}-\frac{2}{3}>\frac{p^2}{6}+\frac{1}{4}
\quad\Rightarrow\quad\Delta=4d-1>\frac{(2p)^2}{6}. 
\end{array}
$$

The case where $p$ is odd ---so we replace $p$ by $2p+1$--- is similar. For $c_1=0$, we consider the shift by $k=-p$, so we obtain: 
$$
\begin{array}{rl}
RR(2p+1)& 
=\binom{\dta + p+1}{2p+1}+\binom{-\dta+p+1}{2p+1}=(-1)^p2\cdot\frac{\ouset{k=0}{p-1}{\prod}(d+k^2)}{(2p)!}
\\[2ex]
\Rightarrow& 
d>\sqrt[2p]{3p}\cdot\frac{4p^2}{\textbf e^2}-\frac{(p-1)(2p-1)}{6}\quad\Rightarrow\quad\Delta=4d>\frac{(2p+1)^2}{6}.
\end{array}
$$
In the remaining case $c_1=1$, shift by $k=-(p+1)$, and we have: 
$$
\begin{array}{rl}
RR(2p+1)& 
=\binom{\dta + p+1/2}{2p+1}+\binom{-\dta+p+1/2}{2p+1}=(-1)^p\cdot\frac{\ouset{k=1}{p}{\prod}(d+k^2-k)}{(2p)!}
\\[2ex]
\Rightarrow& 
d>\sqrt[2p]{3p}\cdot\frac{4p^2}{\textbf e^2}-\frac{p^2-1}{3}\quad\Rightarrow\quad\Delta=4d-1>\frac{(2p+1)^2}{6}.
\end{array}
$$
This concludes the proof. 
\end{m-proof}


\end{document}